\documentclass[12pt]{article}
\usepackage{amssymb, amsmath, amsthm,mathtime}
\usepackage[all]{xy}

\newcommand{\mr}[1]{\mathrm{#1}}
\newcommand{\mf}[1]{\mathfrak{#1}}
\newcommand{\mc}[1]{\mathcal{#1}}

\newcommand{\Z}{{\bf Z}}
\newcommand{\Q}{{\bf Q}}
\newcommand{\C}{{\bf C}}
\newcommand{\zp}{{\bf Z}_p}
\newcommand{\qp}{{\bf Q}_p}
\newcommand{\qbar}{\overline{\Q}}
\newcommand{\qpbar}{\overline{\qp}}

\newcommand{\SL}{\mr{SL}}

\DeclareMathOperator{\Hom}{Hom} 
 \DeclareMathOperator{\Gal}{Gal}

\newtheorem{theorem}{Theorem}
\newtheorem{proposition}[theorem]{Proposition}

\newtheorem{corollary}[theorem]{Corollary}

\newtheorem{conjecture}[theorem]{Conjecture}

\theoremstyle{definition}

\theoremstyle{remark}
\newtheorem*{remark}{Remark}

\renewcommand{\baselinestretch}{1.2}

\begin{document}

\title{Cup products and $L$-values of cusp forms}
\author{Romyar T. Sharifi}
\maketitle

\begin{abstract}
	In this note, we describe a conjecture that, for an odd prime $p$, 
	relates special values of a cup product pairing on cyclotomic $p$-units in 
	$\Q(\mu_p)$ to the $L$-values of newforms satisfying modulo $p$ 
	congruences with Eisenstein series for $\SL_2(\Z)$.
\end{abstract}

\section{Introduction}

Let $p$ denote an odd prime.
In \cite{mcs}, McCallum and the author began a study of a cup product pairing on $p$-units in the $p$th cyclotomic field $F = \Q(\mu_p)$ with values in a Tate twist of the ideal class group of $F$ modulo $p$.  We studied its numerous applications, but the values of this pairing on special cyclotomic $p$-units were quite mysterious to us.  Using the fact that $p$-units of the form $x$ and $1-x$ pair trivially with each other, we were able to narrow down the possibilities for these values for $p$ less than 10,000 (currently, 25,000).  In fact, we were able to find a unique nontrivial possibility for the pairing up to a scalar in the ${\bf F}_p$-group ring of $\Gal(F/\Q)$.  The problem was that the scalar might not be a unit.  Still, we conjectured that it is a unit and that, in general, the pairing had image spanning at least the entire minus part of the class group (see Conjecture \ref{pairconj}).

Very soon after, the author came to realize a connection with the periods of modular forms indirectly through the work of Ihara and Takao \cite{ihara}.  We were led to check the conjectural pairing values on special cyclotomic $p$-units against an old table of Manin's of algebraic periods of cusp forms of small weight \cite{manin}, and we saw that the ratios of the pairing values agreed with the ratios of the periods modulo $p$.  Thanks to an observation of Goncharov, the full correspondence quickly became clear: the values of the cup product pairing correspond in a direct manner to algebraic periods of a newform satisfying a congruence with an Eisenstein series modulo $p$.  This, in turn, led to the author formulate the conjecture (Conjecture \ref{conj}) that is laid out in this article.

Since that time, the author has generalized this conjecture to compare a particular value of a reciprocity map constructed out of an inverse limit of cup products to a two-variable $p$-adic $L$-function of Mazur and Kitagawa.  This and several variants may be found in \cite{me-Lfn}.  Still, we feel that the original conjecture at level 1 and weight $k$ is well worth stating, in that its down-to-earth statement serves as motivation for that work.  For comparison with more recent work, we also state one, cleaner variant (Conjecture \ref{strongconj}).

We remark that recent independent but similar progress has been made on this conjecture by Busuioc \cite{busuioc} and the author \cite[Section 5]{me-Lfn}.  We refer the reader especially to the work of Busuioc in terms of the relation with the conjecture as stated in this article.  She constructs a modular symbol of weight $k$ on $\mr{SL}_2(\Z)$ with $L$-values corresponding to (again, possibly zero) projections of our pairing values.

The author would like thank everyone who gave their advice to him on his way to formulating this conjecture, including, among others, Jordan Ellenberg, Sasha Goncharov, Dick Gross, Yasutaka Ihara, Barry Mazur, Bill McCallum, and William Stein.  He also thanks Celia Busuioc and Glenn Stevens for motivating him to write this article.  This material is based upon work supported by the National Science Foundation under Grant No.\ DMS-0102016 and the Natural Sciences and Engineering Research Council of Canada. 
This research was undertaken, in part, thanks to funding from the Canada Research Chairs Program.

\section{Cup product values} \label{cupproducts}

Let $p$ be an odd prime, and set $F = \Q(\mu_p)$, where $\mu_p$ denotes the group of $p$th roots of unity in an algebraic closure $\qbar$ of $\Q$.  Let $\Delta = \Gal(F/\Q)$.  Let $S$ denote the set consisting of the unique prime above $p$ in $F$, let $G_{F,S}$ denote the Galois group of the maximal unramified outside $p$ extension of $F$, and let $\mc{O}_{F,S}$ denote the ring of $p$-integers in $F$.
McCallum and the author considered the cup product
$$
	H^1(G_{F,S},\mu_p) \otimes H^1(G_{F,S},\mu_p) \xrightarrow{\cup} 
	H^2(G_{F,S},\mu_p^{\otimes 2})
$$
in Galois cohomology. 
(All tensor products in this article are taken over $\Z$.) 
We recall that, by Kummer theory, $H^1(G_{F,S},\mu_p)$
is canonically identified with a subgroup of $F^{\times}/F^{\times p}$ containing
$\mc{O}_{F,S}^{\times}/\mc{O}_{F,S}^{\times p}$.  Furthermore, 
$H^2(G_{F,S},\mu_p^{\otimes 2})$ is canonically identified with
$A_F \otimes \mu_p$ as an ${\bf F}_p[\Delta]$-module, where $A_F$ denotes the $p$-part of the class group of $F$.   We obtain a Galois equivariant 
pairing
$$
	\mc{O}_{F,S}^{\times} \times \mc{O}_{F,S}^{\times} \to A_F \otimes \mu_p.
$$

Now, suppose that $p$ divides the numerator of $B_k/k$ for some positive even $k$, where $B_k$ denotes the $k$th Bernoulli number.  We say that $(p,k)$ is an irregular pair.  Then the eigenspace $A_F^{(1-k)}$ of $A_F$ upon which $\Delta$ acts by the $(1-k)$th power of the Teichm\"uller character is nontrivial.  We have an induced
pairing
$$
	\langle \ \, ,\ \rangle_k \colon 
	\mc{O}_{F,S}^{\times} \times \mc{O}_{F,S}^{\times} \to A_F^{(1-k)} 
	\otimes \mu_p.
$$
McCallum and the author made the following conjecture (see \cite[Conjecture 5.3]{mcs}, but note that Vandiver's conjecture that $A_F^{(q)} = 0$ for all even $q$ is not assumed here):

\begin{conjecture}[McCallum-Sharifi]  \label{pairconj}
	The image of $\langle \ \, ,\ \rangle_k$ spans $A_F^{(1-k)} 
	\otimes \mu_p$.
\end{conjecture}

In particular, the author proved the following \cite[Corollary 5.9]{me-paireis}.

\begin{theorem}
	Conjecture \ref{pairconj} is true for $p < 1000$.
\end{theorem}

Of particular interest are the values of this pairing on cyclotomic $p$-units.
Fix a primitive $p$th root of unity $\zeta_p$.  For positive odd $i$, let
$$
	\eta_i = \prod_{u=1}^{p-1} (1-\zeta_p^u)^{u^{i-1}} \in \mc{O}_{F,S}^{\times}.
$$
For instance, $\eta_1 = p$.
Note that the image of $\eta_i$ modulo $p$th powers lies in the $(1-i)$th eigenspace of 
$\mc{O}_{F,S}^{\times}/\mc{O}_{F,S}^{\times p}$, and it generates it if and only
if $A_F^{(1-i)} = 0$.
Since $\langle \ \, ,\ \rangle_k$ is Galois equivariant, we have
$\langle \eta_i, \eta_j \rangle_k = 0$ if $i + j \not\equiv k \bmod p-1$.  When $i+j \equiv k \bmod p-1$, we set
$$
	e_{i,k} = \langle \eta_i, \eta_j \rangle_k \in A_F^{(1-k)} \otimes \mu_p.
$$
Note that if $j \equiv i \bmod p-1$, then $e_{j,k} = e_{i,k}$.

Under Vandiver's conjecture, which holds for $p <$ 12,000,000 \cite{bcems}, there exists
an isomorphism
$$ 
	\phi \colon A_F^{(1-k)} \otimes \mu_p \xrightarrow{\sim} {\bf F}_p
$$ 
of groups.  Choosing such an isomorphism, we may view the $e_{i,k}$ as elements of ${\bf F}_p$.  For each irregular pair $(p,k)$ with $p <$ 25,000, McCallum and the author have found a unique nonzero possibility for the resulting vector 
$$
	(\phi(e_{1,k}),\phi(e_{3,k}),\ldots,\phi(e_{p-2,k})) \in {\bf F}_p^{(p-1)/2}
$$ 
up to scalar.

\section{The conjecture in terms of $L$-values} \label{Lvalues}

As in Section \ref{cupproducts}, we suppose that $(p,k)$ is an irregular pair, so $k$ is even with $k \ge 12$.  Consider the Eisenstein series of weight $k$ and level $1$:
\begin{eqnarray*}
	G_k = -\frac{B_k}{2k} + \sum_{n=1}^{\infty} \sigma_k(n) q^n,
	&&
	\sigma_k(n) = \sum_{\substack{d \mid n\\d \ge 1}} d^{k-1}.
\end{eqnarray*} 

Let us fix an embedding of $\qbar$ into $\qpbar$, which in particular fixes a 
prime above $p$ in every algebraic extension of $\Q$.
Then there exists newform $f$ of weight $k$ and level $1$ such that
$$
	f \equiv G_k \bmod \mf{p}_f,
$$
where $\mf{p}_f$ is the chosen maximal ideal above $p$ in the ring $\mc{O}_f$ generated over $\Z$ by the coefficients of $f$.  
In fact, it is well-known that the space of weight $k$ cusp forms for $\SL_2(\Z)$ congruent to $G_k$ modulo our fixed prime above $p$ in $\qbar$ has a basis of eigenforms. 

We will use $K_f$ to denote the quotient field of $\mc{O}_f$.
Let us write the $q$-expansion of $f$ as
$$
	f = \sum_{n=1}^{\infty} a_n q^n.
$$
Consider the $L$-values $L(f,j)$ of $f$ for $j$ with $1 \le j \le k-1$, and
set
$$
	\Lambda(f,j) = \frac{(j-1)!}{(-2\pi \sqrt{-1})^j}L(f,j)
$$
for such $j$.  For odd $i$ with $3 \le i \le k-3$, let 
$$
	\rho(f,i) = \Lambda(f,i)/\Lambda(f,1),
$$
which by Eichler-Shimura is an element of $K_f$.

Let $H_f$ be the $\mc{O}_f$-submodule of $K_f$ spanned by the $\rho(f,i)$.
The Eisenstein ideal $I_f$ of $\mc{O}_f$ generated by the $a_n-\sigma_k(n)$ for all $n \ge 1$ is contained in $H_f$, as follows from 
\cite[Theorem 1.3]{manin}.  In fact, the completions of $H_f$ and $I_f$ at $\mf{p}_f$ agree whenever $A_F^{(k-1)} = 0$, as we shall see in Proposition \ref{HI}.

Let $r_{i,f}$ denote the image of $\rho(f,i)$ in $H_f/\mf{p}_fH_f$
for odd $i$ with $3 \le i \le k-3$.
We are now ready to state the first form of our conjecture:

\begin{conjecture} \label{conj}
	Let $\psi \colon H_f/\mf{p}H_f \to {\bf F}_p$.  Then there exists 
	$\phi \colon A_K^{(1-k)} \otimes \mu_p \to {\bf F}_p$ such that
	the two vectors
	$$
		{\bf v}_{\phi} = (\phi(e_{3,k}),\phi(e_{5,k}),\ldots,\phi(e_{k-3,k}))
	$$
	and
	$$
		{\bf w}_{\psi}(f) = (\psi(r_{3,f}),\psi(r_{5,f}),\ldots,\psi(r_{k-3,f}))
	$$
	in ${\bf F}_p^{k/2-2}$ are equal.
\end{conjecture}

Recall from the theory of cyclotomic fields that if $A_F^{(k)} = 0$, then $A_F^{(1-k)}$ is cyclic.  In particular, the hypotheses of the following proposition, which will be proved in Section \ref{complements}, hold under Vandiver's conjecture.

\begin{proposition} \label{onedim}
	If $A_F^{(2-k)} = 0$ and $A_F^{(1-k)}$ is cyclic, then $H_f/\mf{p}_fH_f$
	is one-dimensional over ${\bf F}_p$.  
\end{proposition}

Under the hypotheses of Proposition \ref{onedim}, we then have:

\begin{corollary} \label{vancor}
	If $A_F^{(2-k)} = 0$ and $A_F^{(1-k)}$ is cyclic, then Conjecture 
	\ref{conj} asserts the existence of isomorphisms 
	\begin{eqnarray*}
		\phi \colon A_F^{(1-k)} \otimes \mu_p \xrightarrow{\sim} {\bf F}_p
		& \mr{and} &
		\psi \colon H_f/\mf{p}_f H_f \xrightarrow{\sim} {\bf F}_p
	\end{eqnarray*}
	such that 
	$\phi(e_{i,k}) = \psi(r_{i,f})$
	for all odd $i$ with $3 \le i \le k-3$.
\end{corollary}

Clearly, Corollary \ref{vancor} implies Conjecture \ref{pairconj} under its hypotheses.

\section{The conjecture in terms of modular symbols} \label{modsym}

We continue with our earlier notation.  Let $V_{k-2}$ denote the subgroup of $\Z[X,Y]$ of homogeneous polynomials of degree $k-2$.  This has a left action of $M_2(\Z)^+$ via
$$
	\left( \begin{matrix} a & b \\ c & d \end{matrix} \right) 
	\cdot P(X,Y) = P(dX-bY,-cX+aY)
$$
for $P \in V_{k-2}$.  Let $\Delta$ denote the group of divisors on ${\bf P}^1(\Q)$ and $\Delta_0$ the subgroup of divisors of degree $0$, equipped with the usual action of $M_2(\Z)^+$.  The $\Z$-module $\mc{M}_k(\Z)$ of ``modular symbols'' (see \cite[Chapter 8]{stein}) of weight $k$ on $\SL_2(\Z)$ is then defined as 
$$
	\mc{M}_k(\Z) = (V_{k-2} \otimes \Delta_0)_{\SL_2(\Z)}/\mr{torsion},
$$
i.e., as the maximal torsion-free $\SL_2(\Z)$-invariant quotient of
$V_{k-2} \otimes \Delta_0$.
(We remark here that this is dual to another, perhaps more common definition of the module of $\Z$-valued modular symbols.)  We represent the image of $P \otimes ([b]-[a])$ in this quotient, for $P \in V_{k-2}$ and $a,b \in {\bf P}^1(\Q)$, by $P\{a,b\}$.  In our setting, any modular symbol may be written as an $\Z$-linear combination of symbols of the form $X^{i-1}Y^{k-1-i}\{0,\infty\}$ with $1 \le i \le k-1$. 

For any commutative ring $R$, we set 
$$
	\mc{M}_k(R) = \mc{M}_k(\Z) \otimes R.
$$
If we replace $\Delta_0$ by $\Delta$ in the definition of $\mc{M}_k(R)$, we obtain an $R$-module of boundary symbols $\mc{B}_k(R)$.  We define the space of cuspidal modular symbols $\mc{S}_k(R)$ to be the kernel of the obvious map $\mc{M}_k(R) \to \mc{B}_k(R)$.  
Then $\mc{S}_k(R)$ is generated as an $R$-module by the symbols of the form $X^{i-1}Y^{k-1-i}\{0,\infty\}$ with $2 \le i \le k-2$.  Furthermore, the image of $X^{k-2}\{0,\infty\}$ generates $\mc{B}_k(R)$.
We may also consider the maximal quotients $\mc{M}_k(R)^{\pm}$ and $\mc{S}_k(R)^{\pm}$ upon which the involution 
$$
	P(X,Y)\{0,\infty\} \mapsto P(-X,Y)\{0,\infty\}
$$ 
acts by $\pm 1$.

Let ${\bf T}_k$ denote the cuspidal weight $k$ Hecke algebra for $\SL_2(\Z)$ over $\zp$, and let $\mf{H}_k$ denote the full modular Hecke algebra.  Then $\mf{H}_k$ acts on $\mc{M}_k(\zp)$ and $\mc{B}_k(\zp)$ using the above-described $M_2(\Z)^+$-action, and ${\bf T}_k$ acts on $\mc{S}_k(\zp)$.  
Let $\mc{I}_k$ denote the Eisenstein ideal in ${\bf T}_k$, generated by $T_l-1-l^{k-1}$ for all primes $l$, and let $\mf{m} = p{\bf T}_k + \mc{I}_k$ be the maximal ideal in which it is contained. 

We now present a rather sleeker variant of Conjecture \ref{conj}, which matches better with the discussion in \cite{me-Lfn}.

\begin{conjecture} \label{strongconj}
	There exists an isomorphism 
	\begin{equation*} \label{isom}
		\Xi_k \colon A_F^{(1-k)} \otimes \mu_p \xrightarrow{\sim}
		\mc{S}_k(\zp)^+/\mf{m}\mc{S}_k(\zp)^+
	\end{equation*}
	satisfying 
	$$
	\Xi_k(e_{i,k}) = X^{i-1}Y^{k-i-1}\{0,\infty\}
	$$	
	for odd $i$ with $3 \le i \le k-3$.  
\end{conjecture}

Conjecture \ref{pairconj} is an immediate corollary of Conjecture \ref{strongconj}.

\begin{remark}
	If $A_F^{(2-k)} = 0$, then the work of Kurihara and Harder-Pink shows that the two spaces 
	in Conjecture \ref{strongconj} are indeed isomorphic (see Remark 3.8, Proposition 3.10, and 	
	the footnote on p.\ 292 of \cite{kurihara}). 
\end{remark}

To compare Conjecture \ref{strongconj} with Conjecture \ref{conj}, note that we have a Hecke equivariant integration pairing
\begin{eqnarray*}
	(\ \,,\ )_k \colon \mc{M}_k(\C) \times S_k(\C) \to \C, &&
	(P\{0,\infty\},g)_k = \int_{0}^{\infty} g(z)P(z,1)dz,
\end{eqnarray*}
where the integral is taken along the imaginary axis and $S_k(\C)$ denotes the space of cusp forms of weight $k$ on $\SL_2(\Z)$.  This pairing
is nondegenerate when restricted to $\mc{S}_k(\C)^+$.
Pairing with $f$ induces a period mapping 
$$
	\Phi_f \colon \mc{M}_k(\Z)^+ \to K_f \cdot 
	(X^{k-2}\{0,\infty\},f)_k,
$$  
and we remark that $\Lambda(f,i) = \Phi_f(X^{i-1}Y^{k-1-i}\{0,\infty\})$.

Now, if $\psi \colon H_f/\mf{p}_fH_f \to {\bf F}_p$, then the map
$$
	\psi \circ (\Lambda(f,1)^{-1}\Phi_f) \colon \mc{S}_k(\Z)^+ \to {\bf F}_p
$$ 
extends to a map on $\mc{S}_k(\zp)^+$ that factors through the Eisenstein component.  That is, it induces a map $\psi'_f \colon \mc{S}_k(\zp)_{\mf{m}}^+
\to {\bf F}_p$.  If we define $\phi = \psi'_f \circ \Xi_k$, then we see that
$\phi$ satisfies the statement of Conjecture \ref{conj}.  That is, we have:

\begin{proposition}
	Conjecture \ref{strongconj} implies Conjecture \ref{conj}.
\end{proposition}

The question of whether every ${\bf F}_p$-valued homomorphism on $\mc{S}_k(\zp)^+_{\mf{m}}$ arises as a sum of maps of the form $\psi'_f$ (varying $f$) is
rather more subtle, and we will not attempt to address it here.
\section{Complements} \label{complements}

In this section, we briefly explore the consequences of the work of Kurihara, Harder-Pink, and Ohta for our conjectures. We begin by proving Proposition \ref{onedim}, which claimed that $H_f/\mf{p}_fH_f$ is one-dimensional when 
$A_F^{(2-k)} = 0$ and $A_F^{(1-k)}$ is cyclic.

\begin{proof}[Proof of Proposition \ref{onedim}] 
	Kurihara \cite[Theorem 0.4]{kurihara} has shown that the Eisenstein
	component $({\bf T}_k)_{\mf{m}}$ of ${\bf T}_k$ is Gorenstein. 
	Since
	\begin{equation*} \label{Sk}
		\mc{S}_k(\Z_p)^+_{\mf{m}} \cong \Hom(({\bf T}_k)_{\mf{m}},\zp)
	\end{equation*}   
	(this is standard: see, for instance, \cite[Proposition 2.5]{kurihara}), 
	we have that $\mc{S}_k(\Z_p)^+_{\mf{m}}$ 
	is free of rank $1$ over $({\bf T}_k)_{\mf{m}}$.
	
	Note that $\mc{O}_f/\mf{p}_f$ has residue field ${\bf F}_p$, as the
	completion $({\bf T}_k)_{\mf{m}}$ itself has residue field ${\bf F}_p$ 
	and has the $\mf{p}_f$-completion of $\mc{O}_f$ as a quotient.
	We conclude that $H_f/\mf{p}_f H_f$ is indeed 
	one-dimensional over ${\bf F}_p$, as it is the homomorphic image of
	$\mc{S}_k(\zp)^+_{\mf{m}}$ under the map induced by 
	$\Lambda(f,1)^{-1}\Phi_f$.
\end{proof}

For any $\mc{O}_f$-module $M$, let $M_{\mf{p}_f}$ denote the completion of $M$ at $\mf{p}_f$. Let $\mf{M}$ be the maximal ideal containing the Eisenstein ideal
of $\mf{H}_k$.

\begin{proposition} \label{HI}	
	Suppose that $(\mf{H}_k)_{\mf{M}}$ is Gorenstein.
	Then $(H_f)_{\mf{p}_f} = (I_f)_{\mf{p}_f}.$
\end{proposition}

\begin{proof}
	We need only show that $H_f$ is contained in $(I_f)_{\mf{p}_f}$.
	Consider the usual surjection
	\begin{eqnarray*}
		\Psi_k \colon 
		\mc{M}_k(\zp)^+ \to \zp, && X^{k-2}\{0,\infty\} \mapsto
		1
	\end{eqnarray*}
	taking $\mc{S}_k(\zp)^+$ to zero.  This
	satisfies $\Psi_k(T_l x) = (1+l^{k-1})\Psi_k(x)$ for all primes $l$
	and has kernel $\mc{S}_k(\zp)^+$.
	Since
	$$
		\mc{M}_k(\zp)^+_{\mf{M}} \cong \Hom((\mf{H}_k)_{\mf{M}},\zp),
	$$
	our Gorenstein condition implies that $\mc{M}_k(\zp)^+_{\mf{M}}$ is
	free of rank $1$ over $(\mf{H}_k)_{\mf{M}}$ with generator
	$X^{k-2}\{0,\infty\}$.
	
	Let us define a homomorphism
	$$
		\Phi'_f \colon \mc{M}_k((\mc{O}_f)_{\mf{p}_f})^+ \to (K_f)_{\mf{p}_f}
	$$ 
	by extending $(\mc{O}_f)_{\mf{p}_f}$-linearly the map 
	$\Lambda(f,1)^{-1}\Phi_f$.
	By what we have said above, its image is generated by
	$\Phi'_f(X^{k-2}\{0,\infty\}) = 1$ as an $(\mc{O}_f)_{\mf{p}_f}$-module, 
	so the image is $(\mc{O}_f)_{\mf{p}_f}$.  If we compose
	$\Psi_k$ with the canonical quotient map 
	$$
		\zp \to (\mc{O}_f)_{\mf{p}_f}/(I_f)_{\mf{p}_f},
	$$
	then the resulting map has Hecke eigenvalues equal to the $a_n(f)$ modulo
	$(I_f)_{\mf{p_f}}$.  
	This composition therefore agrees with the composition of $\Phi'_f$ with
	reduction modulo $(I_f)_{\mf{p}_f}$.  
	In particular, since
	$\Psi_k(\mc{S}_k(\zp)^+) = 0$, 
	we have that $H_f$ is contained in $(I_f)_{\mf{p}_f}$.
\end{proof}

We remark that Ohta has shown that $(\mf{H}_k)_{\mf{M}}$ is Gorenstein under
the condition that $A_F^{(k-1)} = 0$, as follows from \cite[Theorem 3.3.2]{ohta}.  One might ask if it is Gorenstein under the weaker condition that $A_F^{(2-k)} = 0$, though we will not explore this further here.

\renewcommand{\baselinestretch}{1}

\end{document}